# ON THE CONTRACTION METHOD WITH DEGENERATE LIMIT EQUATION

By Ralph Neininger[1] and Ludger Rüschendorf

*J. W. Goethe University and Universität Freiburg*

A class of random recursive sequences $(Y_n)$ with slowly varying variances as arising for parameters of random trees or recursive algorithms leads after normalizations to degenerate limit equations of the form $X \stackrel{\mathcal{L}}{=} X$. For nondegenerate limit equations the contraction method is a main tool to establish convergence of the scaled sequence to the "unique" solution of the limit equation. In this paper we develop an extension of the contraction method which allows us to derive limit theorems for parameters of algorithms and data structures with degenerate limit equation. In particular, we establish some new tools and a general convergence scheme, which transfers information on mean and variance into a central limit law (with normal limit). We also obtain a convergence rate result. For the proof we use selfdecomposability properties of the limit normal distribution which allow us to mimic the recursive sequence by an accompanying sequence in normal variables.

**1. Introduction and degenerate limit equations.** A large number of parameters of recursive combinatorial structures, random trees and recursive algorithms satisfy recurrences of the divide-and-conquer type

$$(1) \qquad Y_n \stackrel{\mathcal{L}}{=} \sum_{r=1}^{K} Y_{I_r^{(n)}}^{(r)} + b_n,$$

where $I_r^{(n)}$ are random subgroup sizes in $\{0, \dots, n\}$, $b_n$ is a toll function and $(Y_n^{(r)})_{n \geq 0}$, $r = 1, \dots, K$, are independent copies of the parameter, corresponding to the contribution of subgroup $r$, $\stackrel{\mathcal{L}}{=}$ denotes equality in distribution. Typical parameters $Y_n$ range from the depths and path lengths

---

Received June 2002; revised October 2003.

[1]Supported in part by an Emmy Noether fellowship of the DFG.

of trees, the number of various substructures in combinatorial structures, the number of comparisons, space requirements and other cost measures of algorithms to parameters of communication models, and many more.

The contraction method is an efficient and quite universal probabilistic tool for the asymptotic analysis of recurrences as in (1). It has been introduced for the analysis of the Quicksort algorithm in Rösler (1991) and further developed independently in Rösler (1992, 2001) and Rachev and Rüschendorf (1995), see also the survey article of Rösler and Rüschendorf (2001). It has been applied since then successfully to a large number of problems.

Recently, a fairly general unifying limit theorem for this type of recurrence has been obtained by the contraction method in Neininger and Rüschendorf (2004) in the *nondegenerate case*, where the limit distribution of the normalized recurrence is uniquely characterized by a fixed point equation; we give an illustrative example below. By this result one, in general, obtains the limit distribution from the limiting recurrence and asymptotics of moments.

The aim of this paper is to extend the contraction method and to state a general limit theorem for the *degenerate case*. In the degenerate case the characterizing equations for the normalized algorithm degenerate in the limit to the trivial equation $X \stackrel{\mathcal{L}}{=} X$ and, thus, give no indication on the limit distribution. This case is also quite common in many examples. To simplify the discussion we consider in the first part of the paper recursive sequences $(Y_n)_{n \geq 0}$ which satisfy the distributional recurrence in (1) in the most basic setting, where $K = 1$, that is, we assume that

$$(2) \qquad Y_n \stackrel{\mathcal{L}}{=} Y_{I_n} + b_n, \qquad n \geq n_0,$$

where $n_0 \geq 1$, $(I_n, b_n)$, $(Y_k)$ are independent, $b_n$ is random and $I_n$ is a random index in $\{0, \ldots, n\}$ with $P(I_n = n) < 1$ for $n \geq n_0$. Later on in Section 5 we come back to the more general case as in (1).

To derive a limit in distribution for $(Y_n)$ as in (2) by the contraction method the first step is to introduce a scaling of $Y_n$, say $X_n := (Y_n - \mu_n)/\sigma_n$, where $\mu_n = \mathbb{E}Y_n$ and $\sigma_n = \sqrt{\mathrm{Var}(Y_n)}$ and to derive a recurrence relation for $X_n$:

$$(3) \qquad X_n \stackrel{\mathcal{L}}{=} \frac{\sigma_{I_n}}{\sigma_n} X_{I_n} + b^{(n)}, \qquad n \geq n_0,$$

where

$$b^{(n)} := \frac{1}{\sigma_n}(b_n - \mu_n + \mu_{I_n})$$

and with independence relations as in (2).

The next step to prove a limit theorem for $X_n$ is to establish convergence of the random coefficients in the recursive equation (3):

$$(4) \qquad \frac{\sigma_{I_n}}{\sigma_n} \to A, \qquad b^{(n)} \to b,$$

thus, leading to a limit equation of the form

$$(5) \qquad X \overset{\mathcal{L}}{=} AX + b.$$

Here, $(A, b)$ and $X$ are independent. Essential for the application of the contraction method is that the limit equation (5) has a unique solution under appropriate constraints. The final step of the method is to establish convergence of the $X_n$ to the solution of the limit equation (5).

Many examples of such an approach in the field of analysis of recursive algorithms can be found in Cramer and Rüschendorf (1996), Neininger and Rüschendorf (2004), Rösler (1991, 2001) and Rösler and Rüschendorf (2001) and the references therein.

As a typical example of this approach, we consider the Quickselect algorithm which is designed similarly to the Quicksort algorithm and, as a result, yields a fixed order statistic $x_{(k)}$ of an $n$-tuple of real numbers $x_1, \ldots, x_n$. If $Y_n$ denotes the number of comparisons this algorithm needs to find $x_{(1)}$, then, under the assumption that all permutations of $(x_i)$ are equally likely, $Y_n$ satisfies (2), where $I_n \sim \text{unif}\{0, \ldots, n-1\}$, $b_n = n - 1$, $n_0 = 2$, and $Y_0 = Y_1 = 0$. It is known for this recursion that expectation and variance are of the orders $\mathbb{E}Y_n = 2n + O(1)$ and $\text{Var}(Y_n) = n^2/2 + o(n^2)$, so that, noting that $I_n/n$ has a continuous unif$[0, 1]$ distributed random variable $U$ as its limit, we obtain, after scaling and deriving the limits in (4), a limit equation (5) with $A = U$ and $b = \sqrt{2}(2U - 1)$, thus,

$$(6) \qquad X \overset{\mathcal{L}}{=} UX + \sqrt{2}(2U - 1).$$

The solution of this equation, rescaled by $W = \sqrt{1/2}\, X + 1$, satisfies the equation

$$(7) \qquad W \overset{\mathcal{L}}{=} UW + U,$$

whose unique solution is the Dickman distribution, which is quite common in the analysis of algorithms, as well as in analytic number theory where it originated [see Hwang and Tsai (2002)]. Standard application of the contraction method implies that the fixed point equation (6) has a unique solution $\mathcal{L}(X)$ and that the rescaled quantity $(Y_n - \mathbb{E}Y_n)/\sqrt{\text{Var}(Y_n)}$ converges in distribution to this fixed point.

In this paper we discuss a case which appears quite often for parameters $X_n$ with logarithmic orders for the variance; see the examples below. Here, in the limiting equation (5) we are led to the case $A = 1$, $b = 0$, that is, to the *degenerate limit equation*

$$X \overset{\mathcal{L}}{=} X.$$

The degenerate limit equation does not give any hint to a limit of the recursive sequence $(X_n)$ and so the contraction method does not work in this case.

We will focus in this paper on recursions of the form (2) and the extensions in (1) which lead to a degenerate limit equation and exhibit an asymptotically normal behavior for the scaled quantities $X_n$. We will explain how the normal distribution comes up although the degenerate limit equation does not give any indication for asymptotic normality, and obtain general theorems which lead on the basis of information on mean and variance of $Y_n$ to a central limit law including a rate of convergence. Special cases of our setting are suitable to rederive and extend various limit laws from the field of analysis of algorithms including rates of convergence.

First of all, note that if for $Y_n$ given in (2) we have that $\sigma_n^2 = \text{Var}(Y_n) \sim L(n)$ for $n \to \infty$, with a function $L$ being slowly varying at $\infty$, we obtain

$$\frac{\sigma_{I_n}}{\sigma_n} \sim \sqrt{\frac{L(I_n)}{L(n)}} \to A_1 = 1, \qquad n \to \infty,$$

almost surely, if $I_n$ satisfies mild conditions [see (9)] typically satisfied for applications from the analysis of algorithms. If, furthermore, $b_n$ is appropriately small and $b^{(n)} = \frac{1}{\sigma_n}(b_n - \mu_n + \mu_{I_n}) \to 0$ almost surely, then we are led to the degenerate limit equation for the normalized sequence $(X_n)$. Therefore, degenerate limit equations can be expected for quite general types of recursions.

As an example for the degenerate case, consider the cost $Y_n$ of an unsuccessful search in a random binary search tree as discussed in Cramer and Rüschendorf (1996) and Mahmoud (1992). Here, $(Y_n)$ satisfies (2) with $I_n \sim \text{unif}\{1, \ldots, n - 1\}$, $b_n = 1$, $n_0 = 2$ and $Y_0 = Y_1 = 0$. The toll $b_n = 1$ is small compared to the similar case of Quickselect considered above. In this case the expectation and variance satisfy $\mathbb{E}Y_n = 2 \ln n + O(1)$ and $\text{Var}(Y_n) = 2 \ln n + O(1)$. So the scaling now yields $A = 1$ and $b = 0$ and thus leads to the degenerate limit equation. We come back to this example in Section 4.

Since the case where the variance is a slowly varying function $L(n)$ of the order $(\ln n)^\alpha$ with some $\alpha > 0$ (up to multiplicative constants) is common in the field of analysis of algorithms, we will restrict our setup to this case; for examples see Sections 4 and 5.

The paper is organized as follows: Section 2 contains the basic central limit law, Theorem 2.1. In Section 3 tools are developed to handle degenerate limit equations leading to a proof of Theorem 2.1. In Section 4 as application a couple of limit laws from the field of analysis of algorithms are rederived in a uniform setup. These were previously proven one by one. In the last section we extend our results to obtain central limit theorems for the more complex recurrences of the the divide-and-conquer type in (1). In particular, our limit law covers some more complicated problems related to a maximum-finding algorithm in a broadcast communication model as analyzed in Chen and Hwang (2003).

**2. A central limit law.** Let $(Y_n)_{n \geq 0}$ be a sequence of random variables satisfying the recursion

$$(8) \qquad Y_n \overset{\mathcal{L}}{=} Y_{I_n} + b_n, \qquad n \geq n_0,$$

where $n_0 \geq 1$, $(I_n, b_n), (Y_k)$ are independent, $b_n$ is random and $I_n$ a random index in $\{0, \ldots, n\}$ with $\mathbb{P}(I_n = n) < 1$ for $n \geq n_0$. We denote $\sigma_n = \sqrt{\mathrm{Var}(Y_n)}$ and $\mu_n = \mathbb{E}Y_n$ and use the convention $\ln^\alpha n := (\ln n)^\alpha$ for $\alpha > 0$ and $n \geq 1$. $\|X\|_p$ denotes the $L_p$-norm of a random variable $X$. Then we have the following central limit law, where $\mathcal{N}(0, 1)$ denotes the standard normal distribution.

THEOREM 2.1. *Assume that $(Y_n)_{n \geq 0}$ satisfies the recursion* (8) *with $\|Y_n\|_3 < \infty$ for all $n \geq 0$ and*

$$(9) \qquad \limsup_{n \to \infty} \mathbb{E} \ln \left( \frac{I_n \vee 1}{n} \right) < 0, \qquad \sup_{n \geq 1} \left\| \ln \left( \frac{I_n \vee 1}{n} \right) \right\|_3 < \infty.$$

*Furthermore, assume that for real numbers $\alpha, \lambda, \kappa$ with $0 \leq \lambda < 2\alpha$, the mean and the variance of $Y_n$ satisfy*

$$(10) \qquad \|b_n - \mu_n + \mu_{I_n}\|_3 = O(\ln^\kappa n), \qquad \sigma_n^2 = C \ln^{2\alpha} n + O(\ln^\lambda n),$$

*with some constant $C > 0$. If*

$$(11) \qquad \beta := \tfrac{3}{2} \wedge 3(\alpha - \kappa) \wedge 3(\alpha - \lambda/2) \wedge (\alpha - \kappa + 1) > 1,$$

*then*

$$(12) \qquad \frac{Y_n - \mathbb{E}Y_n}{\sqrt{C} \ln^\alpha n} \overset{\mathcal{L}}{\to} \mathcal{N}(0, 1)$$

*and we have the following rate of convergence for the Zolotarev-metric $\zeta_3$:*

$$(13) \qquad \zeta_3 \left( \frac{Y_n - \mathbb{E}Y_n}{\sqrt{\mathrm{Var}(Y_n)}}, \mathcal{N}(0, 1) \right) = O \left( \frac{1}{\ln^{\beta - 1} n} \right).$$

The Zolotarev metric $\zeta_3$ is defined for distributions $\mathcal{L}(V), \mathcal{L}(W)$ by

$$\zeta_3(\mathcal{L}(V), \mathcal{L}(W)) := \sup_{f \in \mathcal{F}_3} |\mathbb{E}f(V) - \mathbb{E}f(W)|,$$

where $\mathcal{F}_3 := \{f \in C^2(\mathbb{R}, \mathbb{R}) : |f''(x) - f''(y)| \leq |x - y|\}$ is the space of all twice differentiable functions with second derivative being Lipschitz continuous with Lipschitz constant 1. We will use the short notation $\zeta_3(V, W) := \zeta_3(\mathcal{L}(V), \mathcal{L}(W))$. It is well known that convergence in $\zeta_3$ implies weak convergence and that $\zeta_3(V, W) < \infty$ if $\mathbb{E}V = \mathbb{E}W$, $\mathbb{E}V^2 = \mathbb{E}W^2$, and $\|V\|_3, \|W\|_3 < \infty$. The metric $\zeta_3$ is $(3, +)$ ideal, that is, we have for $T$ independent of $(V, W)$ and $c \neq 0$

$$(14) \qquad \zeta_3(V + T, W + T) \leq \zeta_3(V, W), \qquad \zeta_3(cV, cW) = |c|^3 \zeta_3(V, W).$$

For general reference and properties of $\zeta_3$ we refer to Zolotarev (1976, 1977) and Rachev (1991). For implications and interpretation of rates of convergence in the $\zeta_3$ metric see Neininger and Rüschendorf (2002).

**3. Proof of the limit law.** For the scaling of the $Y_n$ we have $\operatorname{Var}(Y_n) \sim C \ln^{2\alpha} n$ with some $\alpha > 0$. Since the scaling of the recurrence requires a scaling for $n = 0, 1$ as well, we define for integers $n \geq 0$ and real $\delta > 0$,

$$L_\delta(n) := \ln(n \vee 1) + \delta \mathbf{1}_{\{0,1\}}(n),$$

where $\mathbf{1}_F$ denotes the indicator function of a set $F$. We use the convention $L_\delta^\alpha(n) := (L_\delta(n))^\alpha$ for $\alpha > 0$.

To prepare for the proof of Theorem 2.1 we provide two calculus lemmas:

LEMMA 3.1.  Let $I_n$ be a random variable in $\{0, \ldots, n\}$ with $\mathbb{P}(I_n = n) < 1$ for all $n$ sufficiently large and with $\limsup_{n \to \infty} \mathbb{E}\ln((I_n \vee 1)/n) < -\varepsilon$ for some $\varepsilon > 0$. Let $(d_n)_{n \geq 0}$, $(r_n)_{n \geq n_0}$ be sequences of nonnegative numbers with

$$d_n \leq \mathbb{E}\left[\left(\frac{L_\delta(I_n)}{L_\delta(n)}\right)^\gamma d_{I_n}\right] + r_n, \qquad n \geq n_0 \geq 2,$$

for some $\gamma > 0$. Then for all $1 < \beta < 1 + \gamma$ and $\delta > 0$ sufficiently small, we have

$$r_n = O\left(\frac{1}{\ln^\beta n}\right) \quad \Longrightarrow \quad d_n = O\left(\frac{1}{\ln^{\beta-1} n}\right).$$

PROOF.  We abbreviate $\eta := \gamma + 1 - \beta$ and choose $\delta = \varepsilon(\eta \wedge 1)/(6\eta)$. There exists an $n_1 \geq n_0$ and an $M > 0$ with $\mathbb{E}\ln((I_n \vee 1)/n) < -\varepsilon$, $p_n := \mathbb{P}(I_n = n) < 1$, $r_n \leq M/\ln^\beta n$, and $(1 + \delta/\ln n)^\eta \leq 1 + 2\eta\delta/\ln n$ for all $n \geq n_1$. We define

$$R := \frac{2M}{\varepsilon(\eta \wedge 1)} \vee \max\{d_k L_\delta^{\beta-1}(k) : 0 \leq k \leq n_1\}$$

and prove $d_n \leq R/L_\delta^{\beta-1}(n)$ by induction. For $0 \leq n \leq n_1$, there is, by definition of $R$, nothing to prove. For $n \geq n_1$, we obtain, using the induction hypothesis,

$$d_n \leq p_n d_n + \mathbb{E}\left[\mathbf{1}_{\{I_n \leq n-1\}}\left(\frac{L_\delta(I_n)}{L_\delta(n)}\right)^\gamma \frac{R}{L_\delta^{\beta-1}(I_n)}\right] + \frac{M}{\ln^\beta n}.$$

This implies

$$(15) \quad d_n \leq \frac{1}{(1-p_n)\ln^{\beta-1} n}\left(R\left(\mathbb{E}\left(\frac{L_\delta(I_n)}{L_\delta(n)}\right)^\eta - p_n\right) + \frac{M}{\ln n}\right)$$

$$(16) \quad \leq \frac{1}{(1-p_n)\ln^{\beta-1} n}\left(R\left(\mathbb{E}\left(1 + \frac{\ln((I_n \vee 1)/n) + \delta}{\ln n}\right)^\eta - p_n\right) + \frac{M}{\ln n}\right).$$

For the estimate of the latter expectation we abbreviate $Z := \ln((I_n \vee 1)/n)$ and the set $F := \{Z > -\delta\}$. Then we have, using $(1-x)^a \leq 1 - ax$ for $x > 0$, $0 < a \leq 1$,

$$\mathbb{E}\left(1 + \frac{Z+\delta}{\ln n}\right)^\eta \leq \mathbb{E}\left[\mathbf{1}_F\left(1 + \frac{\delta}{\ln n}\right)^\eta + \mathbf{1}_{F^c}\left(1 + \frac{Z+\delta}{\ln n}\right)^{\eta \wedge 1}\right]$$

$$\leq \mathbb{E}\left[\mathbf{1}_F\left(1 + \frac{2\eta\delta}{\ln n}\right) + \mathbf{1}_{F^c}\left(1 + \frac{(\eta \wedge 1)(Z+\delta)}{\ln n}\right)\right]$$

$$\leq 1 + \frac{2\eta\delta}{\ln n} + \frac{(\eta \wedge 1)(\mathbb{E}Z + \delta)}{\ln n}.$$

With $\mathbb{E}Z \leq -\varepsilon$ and noting that $\delta \leq \varepsilon(\eta \wedge 1)/(2(2\eta + (\eta \wedge 1)))$, we obtain the estimate

$$\mathbb{E}\left(1 + \frac{Z+\delta}{\ln n}\right)^\eta \leq 1 - \frac{(\eta \wedge 1)\varepsilon}{2\ln n}.$$

Plugging this into (15), we obtain

$$d_n \leq \frac{1}{(1-p_n)\ln^{\beta-1} n}\left(R\left(1 - p_n - \frac{(1 \wedge \eta)\varepsilon}{2\ln n}\right) + \frac{M}{\ln n}\right)$$

$$= \frac{R}{\ln^{\beta-1} n} - \frac{1}{(1-p_n)\ln^\beta n}(R\varepsilon(\eta \wedge 1)/2 - M)$$

$$\leq \frac{R}{\ln^{\beta-1} n},$$

by definition of $R$. $\quad\square$

LEMMA 3.2. *For all $\alpha > 0$ and integers $n \geq 3$ and $1 \leq i \leq n$, we have*

$$\left|\left(\frac{\ln i}{\ln n}\right)^\alpha - 1\right| \leq \frac{2 \vee \alpha}{\ln n}\left|\ln\left(\frac{i}{n}\right)\right|.$$

PROOF.  For $i = 1$, the assertion is true. For $i \geq 2$ and $\alpha \geq 1$, we have, by the mean value theorem, for appropriate $s \in [\ln 2, \ln n]$,

$$\frac{1}{\ln^\alpha n} |\ln^\alpha i - \ln^\alpha n| = \frac{1}{\ln^\alpha n} \alpha s^{\alpha-1} |\ln i - \ln n| \leq \frac{\alpha}{\ln n} \left| \ln \left( \frac{i}{n} \right) \right|.$$

We have

$$\frac{1}{\ln^\alpha n} |\ln^\alpha i - \ln^\alpha n| = \frac{1}{\ln^\alpha n} \frac{|\ln^{2\alpha} i - \ln^{2\alpha} n|}{\ln^\alpha i + \ln^\alpha n} \leq \frac{1}{\ln^{2\alpha} n} |\ln^{2\alpha} i - \ln^{2\alpha} n|.$$

Thus, for $0 < \alpha < 1$, doubling of the exponent $\alpha$ successively yields

$$\frac{1}{\ln^\alpha n} |\ln^\alpha i - \ln^\alpha n| \leq \frac{1}{\ln^{\alpha'} n} |\ln^{\alpha'} i - \ln^{\alpha'} n|$$

with $\alpha' \in [1, 2)$. Then applying the first part implies the assertion.  □

PROOF OF THEOREM 2.1.  We have $\mathbb{E} \ln((I_n \vee 1)/n) < -\varepsilon$ for all $n \geq n_1 \geq n_0$ and some $\varepsilon > 0$. We define the scaled quantities

$$Z_n := \frac{Y_n - \mathbb{E} Y_n}{\sqrt{C} L_\delta^\alpha(n)}, \qquad n \geq 0,$$

with a $\delta > 0$ sufficiently small to be specified later and denote $\tau_n := \sqrt{\mathrm{Var}(Z_n)} = \sigma_n/(\sqrt{C} L_\delta^\alpha(n))$. Thus, we have $\tau_n \to 1$ for $n \to \infty$. The sequence $(Z_n)$ satisfies the recurrence

$$Z_n \overset{\mathcal{L}}{=} \left( \frac{L_\delta(I_n)}{L_\delta(n)} \right)^\alpha Z_{I_n} + b^{(n)}, \qquad n \geq n_1,$$

with

$$b^{(n)} = b^{(n)}(I_n, b_n) = \frac{1}{\sqrt{C} L_\delta^\alpha(n)} (b_n - \mu_n + \mu_{I_n}).$$

Now we define $N_n := \tau_n N$, where $N$ is a standard normal distributed random variable independent of $(I_n, b_n)$, and introduce an accompanying sequence $(Z_n^*)$ by

$$Z_n^* := \left( \frac{L_\delta(I_n)}{L_\delta(n)} \right)^\alpha N_{I_n} + b^{(n)}, \qquad n \geq 0.$$

Note that $Z_n, N_n, Z_n^*$ have identical first and second moment, and finite absolute third moment. Thus, $\zeta_3$ distances between these random variables are finite. We have

(17)                    $\zeta_3(Z_n, N_n) \leq \zeta_3(Z_n, Z_n^*) + \zeta_3(Z_n^*, N_n).$

Using that $\zeta_3$ is $(3, +)$ ideal, compare (14), and conditioning on $(I_n, b_n)$, we obtain

$$
\begin{aligned}
\zeta_3(Z_n, Z_n^*) &= \sup_{f \in \mathcal{F}_3} \left| \int \mathbb{E}\left[ f\left( \left( \frac{L_\delta(k)}{L_\delta(n)} \right)^\alpha Z_k + b^{(n)}(k, s) \right) \right. \right. \\
&\qquad\qquad\qquad \left. \left. - f\left( \left( \frac{L_\delta(k)}{L_\delta(n)} \right)^\alpha N_k + b^{(n)}(k, s) \right) \right] d\mathbb{P}^{(I_n, b_n)}(k, s) \right| \\
&\le \int \zeta_3\left( \left( \frac{L_\delta(k)}{L_\delta(n)} \right)^\alpha Z_k + b^{(n)}(k, s), \right. \\
&\qquad\qquad\qquad \left. \left( \frac{L_\delta(k)}{L_\delta(n)} \right)^\alpha N_k + b^{(n)}(k, s) \right) d\mathbb{P}^{(I_n, b_n)}(k, s) \\
&\le \sum_{k=0}^{n} \mathbb{P}(I_n = k) \left( \frac{L_\delta(k)}{L_\delta(n)} \right)^{3\alpha} \zeta_3(Z_k, N_k).
\end{aligned}
\tag{18}
$$

We will show below that

$$
\zeta_3(Z_n^*, N_n) = O\left( \frac{1}{\ln^\beta n} \right),
\tag{19}
$$

with $\beta$ given in (11). With this estimate, we obtain from (17) and (18) denoting $d_n := \zeta_3(Z_n, N_n)$ and $r_n = \zeta_3(Z_n^*, N_n)$,

$$
\begin{aligned}
d_n &= \zeta_3(Z_n, N_n) \\
&\le \sum_{k=0}^{n} \left( \mathbb{P}(I_n = k) \left( \frac{L_\delta(k)}{L_\delta(n)} \right)^{3\alpha} d_k \right) + r_n \\
&= \mathbb{E}\left[ \left( \frac{L_\delta(I_n)}{L_\delta(n)} \right)^{3\alpha} d_{I_n} \right] + r_n.
\end{aligned}
$$

Lemma 3.1 then implies $d_n = \zeta_3(Z_n, N_n) = O(1/\ln^{\beta-1} n)$ since, by definition of $\beta$, we have $\beta < 1 + 3\alpha$ and $\delta$ can be chosen appropriately. Moreover, we obtain

$$
\begin{aligned}
\zeta_3\left( \frac{Y_n - \mu_n}{\sigma_n}, \mathcal{N}(0, 1) \right) &= \zeta_3\left( \frac{1}{\tau_n} Z_n, \frac{1}{\tau_n} N_n \right) \\
&= \frac{1}{\tau_n^3} \zeta_3(Z_n, N_n) \\
&= O\left( \frac{1}{\ln^{\beta-1} n} \right),
\end{aligned}
$$

which is (13). Since $\zeta_3$ convergence implies weak convergence, we obtain (12).

It remains to establish the bound (19) for $\zeta_3(Z_n^*, N_n)$: We define

$$
G_n := \left( \frac{L_\delta(I_n)}{L_\delta(n)} \right)^\alpha \tau_{I_n},
$$

thus, we have the representation $Z_n^* = G_n N + b^{(n)}$. From $\mathrm{Var}(Z_n^*) = \tau_n^2$, we obtain, in particular, the relation

$$(20) \qquad \mathbb{E}[G_n^2 + (b^{(n)})^2] - \tau_n^2 = 0.$$

Using the closure of the normal familiy under convolution, we have, with the set $A := \{G_n > \tau_n\}$ and its complement $A^c$, the decompositions

$$(21) \qquad Z_n^* \stackrel{\mathcal{L}}{=} \mathbf{1}_A(\tau_n N + \sqrt{G_n^2 - \tau_n^2} N' + b^{(n)}) + \mathbf{1}_{A^c}(G_n N + b^{(n)}),$$

$$(22) \qquad N_n \stackrel{\mathcal{L}}{=} \mathbf{1}_A(\tau_n N) + \mathbf{1}_{A^c}(G_n N + \sqrt{\tau_n^2 - G_n^2} N'),$$

where $N, N', G_n$ are independent, $N \stackrel{\mathcal{L}}{=} N'$. Subsequently, we abbreviate $\Delta_n := |G_n^2 - \tau_n^2|^{1/2}$ and the right-hand sides in (21) and (22) by $\hat{Z}_n^*$ and $\hat{N}_n$, respectively.

We have to estimate $|\mathbb{E}[f(Z_n^*) - f(N_n)]|$ uniformly for $f \in \mathcal{F}_3$. Taylor expansion around $N$ yields $f(x) = f(N) + f'(N)(x - N) + (1/2)f''(N)(x - N)^2 + R(x, N)$ for $x \in \mathbb{R}$. Here we have $|R(x, N)| \leq (1/6)|x - N|^3$ since $f''$ has Lipschitz constant 1. We may subsequently assume that $f''(0) = 0$. If $f''(0) \neq 0$, consider $g(x) := f(x) - (f''(0)/2)x^2$. Then we have $g''(0) = 0$ and, since $Z_n^*, N_n$ have identical second moment, $\mathbb{E}[f(Z_n^*) - f(N_n)] = \mathbb{E}[g(Z_n^*) - g(N_n)]$.

Using the Taylor expansion and representations (21) and (22), we have

$$\mathbb{E}[f(Z_n^*) - f(N_n)] = \mathbb{E}[S_1 + S_2 + R(\hat{Z}_n^*, N) - R(\hat{N}_n, N)],$$

where, for $S_1$, we collect the terms involving the factor $f'(N)$ and, for $S_2$, we collect the terms involving the factor $f''(N)$. Hence, after simplification and using that $N$, $N'$, and $(G_n, b^{(n)})$ are independent, we obtain

$$S_1 = f'(N)(\Delta_n(\mathbf{1}_A - \mathbf{1}_{A^c})N' + b^{(n)}),$$

$$S_2 = \frac{f''(N)}{2}(\Delta_n^2(\mathbf{1}_A - \mathbf{1}_{A^c})(N')^2 + (b^{(n)})^2$$

$$\qquad\qquad + 2b^{(n)}N(\mathbf{1}_A(\tau_n - 1) + \mathbf{1}_{A^c}(G_n - 1))).$$

Since $\mathbb{E}N' = \mathbb{E}b^{(n)} = 0$ and by the independence between $N$ and $b^{(n)}$ and between $N'$ and $(N, G_n)$, we obtain $\mathbb{E}S_1 = 0$. For the estimate of $\mathbb{E}S_2$, first note that we have $\Delta_n^2(\mathbf{1}_A - \mathbf{1}_{A^c}) = G_n^2 - \tau_n^2$. Hence, with (20), the independence of $N, N', (G_n, b^{(n)})$, and $\mathbb{E}(N')^2 = 1$, we obtain

$$\mathbb{E}\frac{f''(N)}{2}(\Delta_n^2(\mathbf{1}_A - \mathbf{1}_{A^c})(N')^2 + (b^{(n)})^2) = 0.$$

Furthermore, note that for $f \in \mathcal{F}_3$ with $f''(0) = 0$, we have

$$|\mathbb{E}[f''(N)N]| = |\mathbb{E}[(f''(N) - f''(0))N]|$$

$$\qquad\qquad \leq \mathbb{E}[|f''(N) - f''(0)||N|] \leq \mathbb{E}N^2 = 1.$$

Thus, with the independence of $N$ to the other quantities, we obtain the bound

$$|\mathbb{E}S_2| \leq \mathbb{E}|b^{(n)}(|\tau_n - 1| + |G_n - 1|)| \leq \|b^{(n)}\|_2(|\tau_n - 1| + \|G_n - 1\|_2).$$

For the remainder terms we have the $O$-estimate

$$\mathbb{E}|R(\hat{Z}_n^*, N)| \leq \tfrac{1}{6}\mathbb{E}[|(\tau_n - 1)N + \Delta_n N' + b^{(n)}|^3 + |(G_n - 1)N + b^{(n)}|^3]$$
$$= O(|\tau_n - 1|^3 + \|\Delta_n\|_3^3 + \|b^{(n)}\|_3^3 + \|G_n - 1\|_3^3).$$

The term $\mathbb{E}|R(\hat{N}_n, N)|$ is bounded by the same $O$-term. Hence, altogether we obtain

$$(23) \qquad \begin{aligned} \zeta_3(Z_n^*, N_n) = {}& O(|\tau_n - 1|^3 + \|\Delta_n\|_3^3 + \|b^{(n)}\|_3^3 + \|G_n - 1\|_3^3 \\ & + \|b^{(n)}\|_2(|\tau_n - 1| + \|G_n - 1\|_2)). \end{aligned}$$

For the estimate of the latter norms and distances note that, using Lemma 3.2, we have

$$\begin{aligned} |G_n - 1| &= \frac{1}{\sqrt{C}\ln^\alpha n}|\sigma_{I_n} - \sqrt{C}\ln^\alpha n| \\ &\leq \frac{1}{C\ln^{2\alpha} n}|\sigma_{I_n}^2 - C\ln^{2\alpha} n| \\ &= \frac{1}{C\ln^{2\alpha} n}|C\ln^{2\alpha}(I_n \vee 1) - C\ln^{2\alpha} n + O(\ln^\lambda n)| \\ &= \left|\ln\left(\frac{I_n \vee 1}{n}\right)\right| O\left(\frac{1}{\ln^{1 \wedge (2\alpha - \lambda)} n}\right). \end{aligned}$$

Analogously, we obtain $|\tau_n - 1| = O(1/\ln^{2\alpha - \lambda} n)$.

With $\sup_{n \geq 1} \|\ln((I_n \vee 1)/n)\|_3 < \infty$, we obtain

$$\|G_n - 1\|_3 = O\left(\frac{1}{\ln^{1 \wedge (2\alpha - \lambda)} n}\right).$$

By definition of $b^{(n)}$ and (10), we have $\|b^{(n)}\|_3 = O(1/\ln^{\alpha - \kappa} n)$. For $\Delta_n$, we obtain

$$\begin{aligned} \|\Delta_n\|_3 &= \|\sqrt{|\tau_n^2 - G_n^2|}\|_3 = \|\tau_n^2 - G_n^2\|_{3/2}^{1/2} \\ &\leq \|\tau_n^2 - G_n^2\|_3^{1/2} \\ &\leq (|\tau_n^2 - 1| + \|G_n^2 - 1\|_3)^{1/2} \\ &= O\left(\frac{1}{\ln^{(1/2) \wedge (\alpha - \lambda/2)} n}\right). \end{aligned}$$

Collecting the estimates, we bound the right-hand side in (23). Estimating there the $L_2$-norms, by $L_3$-norms, we finally obtain

$$\zeta_3(Z_n^*, N_n) = O\left(\frac{1}{\ln^\beta n}\right)$$

with

$$\beta = \tfrac{3}{2} \wedge 3(\alpha - \kappa) \wedge 3(\alpha - \lambda/2) \wedge (\alpha - \kappa + 1) \wedge (3\alpha - \kappa - \lambda).$$

Note that this coincides with the representation for $\beta$ in (11) since we have $3(\alpha - \kappa) \wedge 3(\alpha - \lambda/2) \leq 3\alpha - \kappa - \lambda$. This is seen by distinguishing the cases $\kappa \geq \lambda/2$ and $\kappa < \lambda/2$.  □

In the proof of Theorem 2.1 the limit normal distribution is no longer obtained from the limit fixed-point equation as in the usual contraction method. Instead, as a substitute, the closure of the normal family under convolution used in (21) and (22) allows us to mimic the recurrence satisfied by $(Z_n)$, respectively, by the accompanying sequence $(Z_n^*)$ in terms of normal quantities. This decomposition allows for estimating $\zeta_3(Z_n^*, N_n)$ sufficiently tight. It is easy to see that the scaling property in (21) and (22) essentially characterizes the normal distribution. More precisely, the following lemma explains the occurrence of the normal limit distribution:

LEMMA 3.3 (Characterization of normal distributions).   *Let $X, W$ be independent with mean 0 and variance 1 and assume that for all $q \in (0,1)$,*

(24)
$$X \stackrel{\mathcal{L}}{=} qX + \sqrt{1 - q^2}\, W.$$

*Then we have $X \stackrel{\mathcal{L}}{=} \mathcal{N}(0,1)$.*

PROOF.   From (24) we obtain for all fixed $n \geq 1$, by induction on $1 \leq k \leq n$, that

$$X \stackrel{\mathcal{L}}{=} \sqrt{\frac{n-k+1}{n+1}}X + \sqrt{\frac{1}{n+1}}\sum_{j=1}^{k} W_j,$$

where $W_1, \ldots, W_n, X$ are independent with $W_j \stackrel{\mathcal{L}}{=} W$ for all $j = 1, \ldots, n$. Thus, with $k = n$ we have

$$X \stackrel{\mathcal{L}}{=} \sqrt{\frac{1}{n+1}}X + \sqrt{\frac{n}{n+1}}\left(\frac{1}{\sqrt{n}}\sum_{j=1}^{n} W_j\right).$$

Therefore, the central limit theorem implies $X \stackrel{\mathcal{L}}{=} \mathcal{N}(0,1)$.  □

Note that a similar scaling property valid for stable distributions, in principle, allows the method of proof of Theorem 2.1 to a stable limit theorem.

**4. Applications and discussion.** In this section we give applications of Theorem 2.1. A couple of limit laws obtained before by different means and involving specific calculations for each case are covered by Theorem 2.1:

**Unsuccessful search.** The cost of an unsuccessful search in a random binary search tree with $n$ nodes, as discussed in Cramer and Rüschendorf (1996) and Mahmoud (1992), satisfies recurrence (1) with $I_n \sim \text{unif}\{1, \ldots, n-1\}$, $b_n = 1$ for $n \geq 2$, and $Y_0 = Y_1 = 0$. We have [see Mahmoud (1992)]

$$\mathbb{E}Y_n = 2\ln n + O(1), \qquad \text{Var}(Y_n) = 2\ln n + O(1)$$

and obtain in the notation of Theorem 2.1,

$$\|b_n - \mu_n + \mu_{I_n}\|_3 = \|2\ln(I_n/n) + O(1)\|_3 = O(1).$$

Thus, the parameters in Theorem 2.1 are $\alpha = 1/2$, $\kappa = \lambda = 0$ and we have $\beta = 3/2$. The technical conditions in (9) are satisfied since $\ln((I_n \vee 2)/n) \to \ln U$ in $L_3$ for a unif$[0, 1]$ random variable $U$. (Use representations $I_n = \lceil (n-1)U \rceil$ and decompose the domain of the resulting integral into the intervals $(i/n, (i+1)/n]$ for $i = 0, \ldots, n-1$.) Theorem 2.1 implies the central limit law with a rate of convergence:

$$(25) \qquad \zeta_3\left(\frac{Y_n - \mathbb{E}Y_n}{\sqrt{\text{Var}(Y_n)}}, \mathcal{N}(0, 1)\right) = O\left(\frac{1}{\sqrt{\ln n}}\right).$$

Note that the $1/\sqrt{\ln n}$ rate of convergence for different metrics was shown previously in Cramer and Rüschendorf (1996) based on calculations involving the particular distribution of $I_n$.

**Depths of nodes.** The depth of a random node in a random binary search tree with $n$ nodes satisfies recurrence (1) with $\mathbb{P}(I_n = 0) = 1/n$ and $\mathbb{P}(I_n = k) = 2k/n^2$ for $1 \leq k \leq n-1$ and $b_n = 1$, where $n \geq 2$ and $Y_0 = -1, Y_1 = 0$. We have [see Mahmoud (1992)]

$$\mathbb{E}Y_n = 2\ln n + O(1), \qquad \text{Var}(Y_n) = 2\ln n + O(1),$$

and obtain in the notation of Theorem 2.1,

$$\|b_n - \mu_n + \mu_{I_n}\|_3 = \|2\ln(I_n/n) + O(1)\|_3 = O(1).$$

Hence, the parameters of Theorem 2.1 are given by $\alpha = 1/2$, $\kappa = \lambda = 0$ and we obtain $\beta = 3/2$. The technical conditions in (9) are satisfied since $\ln((I_n \vee 2)/n) \to \ln\sqrt{U}$ in $L_3$ for a unif$[0, 1]$ random variable $U$ and Theorem 2.1 implies the central limit law with a rate of convergence as in (25).

Mahmoud and Neininger (2003) obtained this rate of convergence via an explicit calculation based on the specific distribution of $I_n$ and showed the optimality of the order $1/\sqrt{\ln n}$, that is, $\zeta_3((Y_n - \mu_n)/\sigma_n, \mathcal{N}(0, 1)) =$ $\Theta(1/\sqrt{\ln n})$. This indicates that our estimates in the proof of Theorem 2.1 are tight. See also Mahmoud and Neininger (2003) for a different distributional recurrence satisfied by $(Y_n)$ which leads to the limit equation $X \overset{\mathcal{L}}{=} BX + (1-B)X'$, where $X, X', B$ are independent with $X, X'$ being identically distributed and $B$ Bernoulli(1/2) distributed. This limit equation similar to $X \overset{\mathcal{L}}{=} X$ is as well satisfied by any distribution, hence, also of degenerate type.

**Broadcast communication.** The time $(Y_n)$ of a maximum-finding algorithm for a broadcast communication model with $n$ processors as analyzed in Chen and Hwang [(2003), Algorithm B] satisfies $Y_0 = Y_1 = 1$ and, for $n \geq 2$, recurrence (1) with $I_n \sim \text{unif}\{0, \dots, n-1\}$ and $b_n$ being the time (= number of rounds) used by a leader election algorithm as discussed in Prodinger (1993) and further analyzed in Fill, Mahmoud and Szpankowski (1996). We have [see Chen and Hwang (2003)] $\mathbb{E}b_n^3 = O(\ln^3 n)$ and

$$\mathbb{E}Y_n = \mu \ln^2 n + O(\ln n), \qquad \text{Var}(Y_n) = \sigma^2 \ln^3 n + O(\ln^2 n),$$

with positive constants $\mu, \sigma$. A direct calculation gives, after cancellations of leading terms,

$$\|b_n - \mu_n + \mu_{I_n}\|_3 = O(\ln n).$$

Thus, we have $\alpha = 3/2$, $\kappa = 1$ and $\lambda = 2$, which gives $\beta = 3/2$. This implies the following

COROLLARY 4.1. *The time* $(Y_n)$ *of Algorithm* B *in Chen and Hwang (2003), as introduced above, satisfies*

$$\zeta_3\left(\frac{Y_n - \mathbb{E}Y_n}{\sqrt{\text{Var}(Y_n)}}, \mathcal{N}(0,1)\right) = O\left(\frac{1}{\sqrt{\ln n}}\right).$$

The same bound for the rate for the Kolmogorov metric was obtained in Chen and Hwang (2003).

**5. Extensions and applications.** We consider now the more general recurrence for $(Y_n)$, as in (1),

$$(26) \qquad Y_n \overset{\mathcal{L}}{=} \sum_{r=1}^{K} Y_{I_r^{(n)}}^{(r)} + b_n, \qquad n \geq n_0,$$

where $n_0, K \geq 1$, $b_n$ is a random variable, $I_1^{(n)}, \dots, I_K^{(n)} \in \{0, \dots, n\}$ are random indices, and $(Y_k^{(1)}), \dots, (Y_k^{(K)})$ distributional copies of $(Y_k)$ such that $(Y_k^{(1)}), \dots, (Y_k^{(K)})$, $(I_1^{(n)}, \dots, I_K^{(n)}, b_n)$ are independent. Many examples of *divide-and-conquer* type algorithms lead to this equation and have been considered in the analysis of algorithms literature.

We introduce the scaling $X_n := (Y_n - \mu_n)/\sigma_n$, where $\mu_n = \mathbb{E}Y_n$ and $\sigma_n = \sqrt{\mathrm{Var}(Y_n)}$ and obtain as recurrence relation for $X_n$,

$$(27) \qquad X_n \overset{\mathcal{L}}{=} \sum_{r=1}^{K} \frac{\sigma_{I_r^{(n)}}}{\sigma_n} X_{I_r^{(n)}}^{(r)} + b^{(n)}, \qquad n \geq n_0,$$

where

$$b^{(n)} := \frac{1}{\sigma_n}\left(b_n - \mu_n + \sum_{r=1}^{K} \mu_{I_r^{(n)}}\right)$$

and $(X_k^{(1)}), \ldots, (X_k^{(K)}), (I_1^{(n)}, \ldots, I_K^{(n)}, b_n)$ are independent, $(X_k^{(1)}), \ldots, (X_k^{(K)})$ being distributional copies of $(X_k)$.

Extensions of Theorem 2.1 in various directions are possible. We give as an example a theorem tailored for the case when the coefficients $\sigma_{I_r^{(n)}}/\sigma_n$ in (27) behave roughly as follows:

$$\frac{\sigma_{I_1^{(n)}}}{\sigma_n} \to A_1 = 1, \qquad \frac{\sigma_{I_r^{(n)}}}{\sigma_n} \to A_r = 0, \qquad r = 2, \ldots, K.$$

We assume that $\limsup_{n \to \infty} \sum_{r=1}^{K} \mathbb{P}(I_r^{(n)} = n) < 1$ and denote $\sigma_n = \sqrt{\mathrm{Var}(Y_n)}$ and $\mu_n = \mathbb{E}Y_n$.

THEOREM 5.1. *Assume that* $(Y_n)_{n \geq 0}$ *satisfies the recurrence* (26) *with* $\|Y_n\|_3 < \infty$ *for all* $n \geq 0$, *and*

$$(28) \quad \limsup_{n \to \infty} \mathbb{E}\ln\left(\frac{1}{n}\prod_{r=1}^{K}(I_r^{(n)} \vee 1)\right) < 0, \qquad \sup_{n \geq 1}\left\|\ln\left(\frac{I_1^{(n)} \vee 1}{n}\right)\right\|_3 < \infty.$$

*Furthermore, assume that for real numbers* $\alpha, \lambda, \kappa$ *with* $0 \leq \lambda < 2\alpha$, *the mean and the variance of* $Y_n$ *satisfy*

$$\left\|b_n - \mu_n + \sum_{r=1}^{K} \mu_{I_r^{(n)}}\right\|_3 = O(\ln^\kappa n), \qquad \sigma_n^2 = C\,\ln^{2\alpha} n + O(\ln^\lambda n),$$

*with some constant* $C > 0$ *and that for some real number* $\xi \geq 0$, *we have*

$$\|\ln^\alpha(I_r^{(n)} \vee 1)\|_3 = O(\ln^\xi n), \qquad r = 2, \ldots, K.$$

*If*

$$(29) \qquad \beta := \tfrac{3}{2} \wedge 3(\alpha - \kappa) \wedge 3(\alpha - \xi) \wedge 3(\alpha - \lambda/2) \wedge (\alpha - \kappa + 1) > 1,$$

*then*

$$(30) \qquad \frac{Y_n - \mathbb{E}Y_n}{\sqrt{C}\,\ln^\alpha n} \overset{\mathcal{L}}{\to} \mathcal{N}(0, 1),$$

*and we have the following rate of convergence for the Zolotarev-metric $\zeta_3$:*

$$\zeta_3\left(\frac{Y_n - \mathbb{E}Y_n}{\sqrt{\text{Var}(Y_n)}}, \mathcal{N}(0,1)\right) = O\left(\frac{1}{\ln^{\beta-1} n}\right).$$

For the proof we need a substitute for Lemma 3.1:

LEMMA 5.2.   *Let* $I_1^{(n)}, \ldots, I_K^{(n)}$ *be random variables in* $\{0, \ldots, n\}$ *with* $\sum_{r=1}^{K} \mathbb{P}(I_r^{(n)} = n) < 1$ *for all $n$ sufficiently large, and*

$$\limsup_{n \to \infty} \mathbb{E}\ln\left(\frac{1}{n}\prod_{r=1}^{K}(I_r^{(n)} \vee 1)\right) < -\varepsilon$$

*for some* $\varepsilon > 0$. *Let* $(d_n)_{n \geq 0}$, $(r_n)_{n \geq n_0}$ *be sequences of nonnegative numbers with*

$$d_n \leq \mathbb{E}\left[\sum_{r=1}^{K}\left(\frac{L_\delta(I_r^{(n)})}{L_\delta(n)}\right)^\gamma d_{I_r^{(n)}}\right] + r_n, \qquad n \geq n_0,$$

*for some* $\gamma > 1$. *Then, for all* $1 < \beta \leq \gamma$ *and* $\delta > 0$ *sufficiently small, we have*

$$r_n = O\left(\frac{1}{\ln^\beta n}\right) \quad \Longrightarrow \quad d_n = O\left(\frac{1}{\ln^{\beta-1} n}\right).$$

The proof of Lemma 5.2 follows the argument of the proof of Lemma 3.1. Note that we have the more restrictive condition $1 \leq \beta \leq \gamma$ compared to $1 \leq \beta \leq \gamma + 1$ in Lemma 3.1. This allows for replacing the analog of the estimates (15) and (16) in the proof of Lemma 3.1 by

$$\mathbb{E}\sum_{r=1}^{K}\left(\frac{\ln(I_r^{(n)} \vee 1)}{\ln n}\right)^\eta = \mathbb{E}\left[\left(1 + \frac{\ln((I_1^{(n)} \vee 1)/n)}{\ln n}\right)^\eta + \sum_{r=2}^{K}\left(\frac{\ln(I_r^{(n)} \vee 1)}{\ln n}\right)^\eta\right]$$

$$\leq \mathbb{E}\left[1 + \eta\frac{\ln((I_1^{(n)} \vee 1)/n)}{\ln n} + \sum_{r=2}^{K}\frac{\ln(I_r^{(n)} \vee 1)}{\ln n}\right]$$

$$\leq 1 + \frac{1}{\ln n}\mathbb{E}\ln\left(\frac{1}{n}\prod_{r=1}^{K}(I_r^{(n)} \vee 1)\right).$$

For this we used that $\eta = \gamma + 1 - \beta \geq 1$.

We sketch the extension of the techniques for Theorem 2.1 to obtain Theorem 5.1.

PROOF OF THEOREM 5.1.   (Sketch) We have $\mathbb{E}\ln((1/n)\prod_{r=1}^{K}(I_r^{(n)} \vee 1)) < -\varepsilon$ for all $n \geq n_1$ with an appropriate $\varepsilon > 0$ and $n_1 \geq n_0$. With $\delta > 0$ sufficiently small, we define the scaled quantities

$$Z_n := \frac{Y_n - \mathbb{E}Y_n}{\sqrt{C}L_\delta^\alpha(n)}, \qquad n \geq 0,$$

and denote $\tau_n := \sqrt{\operatorname{Var}(Z_n)} = \sigma_n/(\sqrt{C}L_\delta^\alpha(n))$; thus, $\tau_n \to 1$ for $n \to \infty$. We have the recurrence

$$Z_n \stackrel{\mathcal{L}}{=} \sum_{r=1}^K \left(\frac{L_\delta(I_r^{(n)})}{L_\delta(n)}\right)^\alpha Z_{I_r^{(n)}}^{(r)} + b^{(n)}, \qquad n \geq n_1,$$

with

$$b^{(n)} = \frac{1}{\sqrt{C}L_\delta^\alpha(n)}\left(b_n - \mu_n + \sum_{r=1}^K \mu_{I_r^{(n)}}\right).$$

We define $N_n^{(r)} := \tau_n N^{(r)}$, where $N^{(1)}, \ldots, N^{(K)}$ are standard normal distributed random variables such that $(I_n, b_n), N^{(1)}, \ldots, N^{(K)}$ are independent. Also, we introduce an accompanying sequence $(Z_n^*)$ by

$$Z_n^* := \sum_{r=1}^K \left(\frac{L_\delta(I_r^{(n)})}{L_\delta(n)}\right)^\alpha N_{I_r^{(n)}}^{(r)} + b^{(n)}, \qquad n \geq n_1.$$

Then, with $d_n = \zeta_3(Z_n, N_n)$ and $r_n = \zeta_3(Z_n^*, N_n^{(1)})$, we obtain similarly to the argument in the proof of Theorem 2.1,

$$d_n \leq \mathbb{E}\left[\sum_{r=1}^K \left(\frac{L_\delta(I_r^{(n)})}{L_\delta(n)}\right)^{3\alpha} d_{I_r^{(n)}}\right] + r_n.$$

For the estimate of $r_n = \zeta_3(Z_n^*, N_n)$, we define

$$G_n := \left(\frac{L_\delta(I_1^{(n)})}{L_\delta(n)}\right)^\alpha \tau_{I_1^{(n)}}, \qquad A := \{G_n > \tau_n\}, \qquad \Delta_n := \sqrt{|G_n^2 - \tau_n^2|}$$

and use the representations

$$Z_n^* \stackrel{\mathcal{L}}{=} \mathbf{1}_A\left(\tau_n N^{(1)} + \Delta_n N' + \sum_{r=2}^K \left(\frac{L_\delta(I_r^{(n)})}{L_\delta(n)}\right)^\alpha N_{I_r^{(n)}}^{(r)} + b^{(n)}\right)$$

$$+ \mathbf{1}_{A^c}\left(G_n N^{(1)} + \sum_{r=2}^K \left(\frac{L_\delta(I_r^{(n)})}{L_\delta(n)}\right)^\alpha N_{I_r^{(n)}}^{(r)} + b^{(n)}\right),$$

$$N_n \stackrel{\mathcal{L}}{=} \mathbf{1}_A(\tau_n N^{(1)}) + \mathbf{1}_{A^c}(G_n N^{(1)} + \Delta_n N'),$$

where $N'$ is standard normal distributed and independent of the other random variates. With corresponding estimates, as in the proof of Theorem 1.1, we find

$$\zeta_3(Z_n^*, N_n) = O\left(|\tau_n - 1|^3 + \|\Delta_n\|_3^3 + \|b^{(n)}\|_3^3 + \|G_n - 1\|_3^3\right.$$

$$+ \|b^{(n)}\|_2(|\tau_n - 1| + \|G_n - 1\|_2) + \frac{1}{\ln^{3\alpha} n} \sum_{r=2}^{K} \|\ln^\alpha(I_r^{(n)} \vee 1)\|_3^3 \Big)$$

$$= O\Big(\frac{1}{\ln^\beta n}\Big),$$

with $\beta$ given in (29). Since $\beta \le 3\alpha$, Lemma 5.2 completes the proof. $\square$

As applications of Theorem 5.1 we discuss various cost measures $(Y_n)$ for a maximum finding algorithm in a broadcast communication model with $n$ processors as analyzed in Chen and Hwang [(2003), Algorithm A]. We use their expansions for mean and variance and, by Theorem 5.1, rederive the central limit laws. Additionally, we endow them with new rates of convergences. Several cost measures $(Y_n)$ of this algorithm satisfy the recurrence

$$(31) \qquad Y_n \stackrel{\mathcal{L}}{=} Y_{I_1^{(n)}}^{(1)} + Y_{I_2^{(n)}}^{(2)} + b_n, \qquad n \ge 2,$$

with relations as in (26), where $b_n$ varies for different cost measures, whereas the distribution of the indices $(I_1^{(n)}, I_2^{(n)})$ is in all cases given by

$$\mathbb{P}((I_1^{(n)}, I_2^{(n)}) = (j, k)) = \begin{cases} 2^{-n}, & (j,k) = (0,0), \\ \binom{n-k-1}{j-1} 2^{-n}, & k \ge 0, 1 \le j \le n-k. \end{cases}$$

In particular, we have that the marginal $I_1^{(n)}$ is binomial $B(n, 1/2)$ distributed and

$$\mathbb{P}(I_2^{(n)} = k) = \begin{cases} \frac{1}{2} + 2^{-n}, & k = 0, \\ 2^{-(k+1)}, & 1 \le k \le n-1. \end{cases}$$

The technical conditions in Theorem 5.1 regarding the indices $(I_1^{(n)}, I_2^{(n)})$ are, hence, satisfied: We have $\mathbb{P}(I_1^{(n)} = n) + \mathbb{P}(I_2^{(n)} = n) = 2^{-n} < 1$ for all $n \ge 1$ and $\|\ln((I_1^{(n)} \vee 1)/n)\|_3 \to \ln 2$ since $I_1^{(n)}$ is binomial $B(n, 1/2)$ distributed, thus, we have $\sup_{n \ge 1} \|\ln((I_1^{(n)} \vee 1)/n)\|_3 < \infty$. For the verification of the first condition in (28) note that we have $\mathbb{E} \ln((I_1^{(n)} \vee 1)/n) \to -\ln 2$ and, therefore, it is sufficient to show $\limsup_{n \to \infty} \mathbb{E} \ln(I_2^{(n)} \vee 1) < \ln 2$. We have

$$\mathbb{E} \ln(I_2^{(n)} \vee 1) = \sum_{k=2}^{n-1} \frac{\ln k}{2^{k+1}} \le \frac{\ln 2}{8} + \sum_{k=3}^{\infty} \frac{k}{2^{k+1}} = \frac{\ln 2}{8} + \frac{1}{2} < 0.6 < \ln 2.$$

Chen and Hwang (2003) analyze three cost measures, namely, the time (= number of rounds) taken by the algorithm, the number of coin flips performed and the number of comparisons performed. The number of coin flips does not lead to a degenerate limit equation $X \stackrel{\mathcal{L}}{=} X$ and can be treated by standard application of the contraction method, see, for example, Rösler (2001). We focus on the other two more delicate parameters:

**Time of the algorithm.** The time $(Y_n)$ of the maximum finding Algorithm A analyzed in Chen and Hwang (2003) satisfies (31) with $b_n = 1$ and $Y_0 = Y_1 = 1$. Mean and variance satisfy, see Chen and Hwang (2003),

$$\mathbb{E}Y_n = \hat{\mu}\ln n + O(1), \qquad \mathrm{Var}(Y_n) = \hat{\sigma}^2\ln n + O(1),$$

with constants $\hat{\mu}, \hat{\sigma} > 0$ being explicitly known. Hence, in the notation of Theorem 5.1 we have

$$\|b_n - \mu_n + \mu_{I_1^{(n)}} + \mu_{I_2^{(n)}}\|_3 = O(1), \qquad \|\ln^{1/2}(I_2^{(n)} \vee 1)\|_3 = O(1),$$

using that also $\|\ln(I_2^{(n)} \vee 1)\|_3 = O(1)$. Thus, we have $\alpha = 1/2$ and $\kappa = \lambda = \xi = 0$, which gives $\beta = 3/2$. With Theorem 5.1 we rederive the central limit law and add the following rate of convergence:

COROLLARY 5.3. *The time (= number of rounds) $(Y_n)$ of the maximum finding Algorithm A in Chen and Hwang (2003), as introduced above, satisfies*

$$\zeta_3\left(\frac{Y_n - \mathbb{E}Y_n}{\sqrt{\mathrm{Var}(Y_n)}}, \mathcal{N}(0,1)\right) = O\left(\frac{1}{\sqrt{\ln n}}\right).$$

**Number of comparisons.** The number of comparisons $(Y_n)$ of the maximum finding Algorithm A was analyzed in Chen and Hwang (2003). It satisfies (31) with $b_n = n - I_1^{(n)}$ and $Y_0 = Y_1 = 0$. Mean and variance have the expansions, see Chen and Hwang (2003),

$$\mathbb{E}Y_n = n + \bar{\mu}\ln n + O(1), \qquad \mathrm{Var}(Y_n) = \bar{\sigma}^2\ln n + O(1),$$

with constants $\bar{\mu}, \bar{\sigma} > 0$ being explicitly known. Hence, in the notation of Theorem 5.1 we obtain, after cancelation,

$$\|b_n - \mu_n + \mu_{I_1^{(n)}} + \mu_{I_2^{(n)}}\|_3 = \|\bar{\mu}\ln((I_1^{(n)} \vee 1)/n) + I_2^{(n)} + \bar{\mu}\ln I_2^{(n)} + O(1)\|_3$$
$$= O(1).$$

Thus, we have $\alpha = 1/2$ and $\kappa = \lambda = \xi = 0$, which gives $\beta = 3/2$. Theorem 5.1 rederives the central limit law and adds a rate of convergence:

COROLLARY 5.4. *The number of comparisons $(Y_n)$ of the maximum finding Algorithm A in Chen and Hwang (2003), as introduced above, satisfies*

$$\zeta_3\left(\frac{Y_n - \mathbb{E}Y_n}{\sqrt{\mathrm{Var}(Y_n)}}, \mathcal{N}(0,1)\right) = O\left(\frac{1}{\sqrt{\ln n}}\right).$$

**Acknowledgment.** The authors thank two anonymous reviewers for several valuable suggestions on this paper.

DEPARTMENT OF MATHEMATICS
J.W. GOETHE UNIVERSITY
ROBERT-MAYER-STR. 10
60325 FRANKFURT A.M.
GERMANY
E-MAIL: neiningr@math.uni-frankfurt.de
URL: www.math.uni-frankfurt.de/neiningr/

INSTITUT FÜR MATHEMATISCHE STOCHASTIK
UNIVERSITÄT FREIBURG
ECKERSTR. 1
79104 FREIBURG
GERMANY
E-MAIL: ruschen@stochastik.uni-freiburg.de
URL: www.stochastik.uni-freiburg.de/homepages/rueschendorf/